\documentclass[12pt]{amsart}
\usepackage{amssymb}
\usepackage{amsfonts}
\usepackage{amsmath}
\usepackage{array}
\usepackage{rotating}
\usepackage{epsf}

\newtheorem{example}{Example}[section]

\newtheorem{theorem}[example]{Theorem}

\newtheorem{corollary}[example]{Corollary}

\newtheorem{lemma}[example]{Lemma}

\def\Proof{\noindent \it Proof -- \rm}
\def\qed{\hspace{3.5mm} \hfill \vbox{\hrule height 3pt depth 2 pt width 2mm}
\bigskip}

\def\dd{\stackrel{\leftarrow}{\partial}}
\def\last{{\rm last}}
\def\des{{\rm des}}
\def\maj{{\rm maj}}
\def\desris{{\rm desris}}

\def\coimaj{{\rm coimaj}}
\def\tmaj{\theta{\rm maj}}
\def\tAdj{\theta{\rm Adj}}
\def\tadj{\theta{\rm adj}}

\def\E{{\mathbb E}}
\def\J{{\bf J}}
\def\FQSym{{\bf FQSym}}

\def\Std{{\rm Std}}

\def\<{\langle}
\def\>{\rangle}

\def\C{\operatorname{\mathbb C}}
\def\Z{\operatorname{\mathbb Z}}
\def\N{\operatorname{\mathbb N}}

\def\F{{\bf F}}

\def\G{{\bf G}}

\def\SG{{\mathfrak S}}

\def\A{{\bf A}}
\def\AA{{\mathcal A}}

\def\Sym{{\bf Sym}}

\def\Des{\operatorname{Des}}

\def\J{{\mathbf J}}
\def\JJ{{\mathbb J}}

\def\shuff#1#2{\mathbin{
\hbox{\vbox{ \hbox{\vrule \hskip#2 \vrule height#1 width 0pt
}%
\hrule}%
\vbox{ \hbox{\vrule \hskip#2 \vrule height#1 width 0pt
\vrule }%
\hrule}%
}}}

\def\biw#1#2{\left[\,\begin{matrix}#1 \cr #2\end{matrix}\,\right]}

\def\empile#1#2{{\scriptsize\scriptstyle \begin{matrix}#1 \cr #2\end{matrix}}}

\def\shuf{{\mathchoice{\shuff{7pt}{3.5pt}}%
{\shuff{6pt}{3pt}}%
{\shuff{4pt}{2pt}}%
{\shuff{3pt}{1.5pt}}}}%
\def\shuffle{\,\shuf\,}

\def\Tabvrule{\vrule width-0.4pt}       
\def\Tabhrule{\hrule \hrule height-0.4pt} 
\def\Tabstrut{\vrule height2.2ex 
                     depth0.8ex  
                     width0ex    
\relax}

\def\PasCase#1{\omit%
            $\vcenter{\hbox {\vbox to 0.4pt{}}
               \hbox{\makebox[3ex]{\Tabstrut$#1$}}}%
               \Tabvrule$}
\def\PasCasePoint{\PasCase{\cdot}}
\def\DessinCarre#1{%
    \vcenter{\hbox{}\hrule
             \hbox{\vrule\makebox[3ex]{\Tabstrut$#1$}\vrule}\Tabhrule}%
             \Tabvrule}
\def\GenRuban#1{\vcenter{\halign{&$\DessinCarre{##}$\cr#1}}\egroup}

\def\sTabvrule{\vrule width-0.4pt}
\def\sTabhrule{\hrule \hrule height-0.4pt}
\def\sTabstrut{\vrule height1.6ex depth0.6ex width0ex \relax}
\def\sDessinCarre#1{%
    \vcenter{\hbox{}\hrule
             \hbox{\vrule\makebox[2.3ex]%
                  {\sTabstrut$\scriptstyle#1$}\vrule}\sTabhrule}%
             \sTabvrule}
\def\sGenRuban#1{\vcenter{\halign{&$\sDessinCarre{##}$\cr#1}}\egroup}

\def\ruban{%
  \bgroup
  \let\ =\omit
  \let\\=\cr
  \let\x=\times
  \let\.=\PasCasePoint
  \offinterlineskip
  \GenRuban}

\def\sruban{%
  \bgroup
  \let\ =\omit
  \let\x=\times
  \let\\=\cr
  \offinterlineskip
  \sGenRuban}

\title{Noncommutative Bessel Symmetric Functions}

\author[J.-C.~Novelli and J.-Y.~Thibon]%
{Jean-Christophe Novelli and Jean-Yves Thibon}

\address[] {Institut Gaspard Monge, Universit\'e de Marne-la-Vall\'ee \\
5 Boulevard Descartes \\Champs-sur-Marne \\77454 Marne-la-Vall\'ee cedex 2 \\
FRANCE}
\email[Jean-Christophe Novelli]{novelli@univ-mlv.fr}
\email[Jean-Yves Thibon]{jyt@univ-mlv.fr} 
\date{}

\begin{document}

\begin{abstract}
The consideration of tensor products of $0$-Hecke algebra modules
leads to natural analogs of the Bessel $J$-functions in the algebra
of noncommutative symmetric functions. This provides a simple explanation
of various combinatorial properties of  Bessel functions.
\end{abstract}

\maketitle


\section{Introduction}

It is known that the theory of noncommutative symmetric functions and
quasi-symmetric functions is related to $0$-Hecke algebras in the same
way as ordinary symmetric functions are related to symmetric groups.
Thus, one may expect that natural questions about representations of
$0$-Hecke algebras lead to the introduction of interesting families
of noncommutative symmetric functions. By ``interesting'', one may
mean ``noncommutative analogs''  of the Frobenius characteristics 
of representations of symmetric groups based on combinatorial objects, 
which may themselve give back various identities for the ordinary,
exponential, or $q$-exponential generating functions of these objects.
This amounts to specialize the complete noncommutative symmetric functions
$S_n(A)$ to $h_n(X)$, $1$, $\frac{1}{n!}$ or $\frac{1}{(q)_n}$, respectively.

Examples of this situation can be found in \cite{NT3}, where the
analysis of the representation of $H_n(0)$ on parking functions
leads naturally to the combinatorics of the noncommutative Lagrange
inversion formula, and to the introduction of noncommutative analogs
of various special functions, such as the Abel polynomials, the Lambert
binomial series or the Eisenstein exponential, and allows one to recover
in a straightforward and unified way a number of enumerative formulas.

The present paper addresses the following question. The $0$-Hecke
algebra is the algebra of a monoid, hence admits a natural coproduct
for which the monoid elements are grouplike. This allows one to define
the tensor product of $0$-Hecke modules, which induces on
quasi-symmetric functions an analog of the
internal product of symmetric functions. What are the properties of this
operation, and of the dual coproduct on noncommutative symmetric functions?

It turns out that  the second part of the question is the most interesting.
Basically, the answer is: {\em the dual coproduct
governs the combinatorics of Bessel functions}. Indeed, its explicitation
leads to the introduction of noncommutative analogs $\J_n(A,B)$ of
the $J$-functions of integer index, of which a few basic
properties are readily established. Then, the above mentioned 
specializations (and other more complicated ones) give back
various classical enumerative formulas.

\section{Background}

\subsection{Notations}
Our notations for noncommutative symmetric functions are as in
\cite{NCSF1,NCSF2}. The Hopf algebra of noncommutative symmetric
functions is denoted by $\Sym$, or by $\Sym(A)$ if we consider the realization
in terms of an auxiliary alphabet. Bases of $\Sym_n$ are labelled by
compositions $I$ of $n$. The noncommutative complete and elementary functions
are denoted by $S_n$ and $\Lambda_n$, and the notation $S^I$ means
$S_{i_1}\cdots S_{i_r}$. The ribbon basis is denoted by $R_I$.
The notation $I\vDash n$ means that $I$ is a composition of $n$.
The conjugate composition is denoted by $I^\sim$.

The graded dual of $\Sym$ is $QSym$ (quasi-symmetric functions).
The dual basis of $(S^I)$ is $(M_I)$ (monomial), and that of $(R_I)$
is $(F_I)$.

The {\em Hecke algebra} $H_n(q)$ ($q\in\C$) is the $\C$-algebra
generated by $n-1$ elements $T_1,\ldots,T_{n-1}$ satisfying
the braid relations and $(T_i-1)(T_i+q)=0$. We are interested in the case
$q=0$, whose representation theory can be described in terms of
quasi-symmetric functions and noncommutative symmetric functions
\cite{NCSF4,NCSF6}.

The Hopf structures on $\Sym$ and $QSym$ allows one 
to extend the $\lambda$-ring notation 
of ordinary symmetric functions (see \cite{NCSF2}, and \cite{Las}
for background on the original commutative version).
If $A$ and $X$  totally ordered sets of noncommuting and commuting variables
respectively, the noncommutative symmetric functions of $XA$  are defined by
\begin{equation}
\sigma_t(XA)=\sum_{n\ge 0}t^n S_n(XA) =
\prod_{x\in X}^{\rightarrow}\sigma_{tx}(A)= \sum_I t^{|I|}M_I(X)S^I(A)\,.
\end{equation}
Thanks to the commutative image
homomorphism $\Sym\rightarrow Sym$,
 noncommutative symmetric functions can be evaluated
on any element $x$ of a $\lambda$-ring, $S_n(x)$ being $S^n(x)$, the $n$-th
symmetric power. Recall that $x$ is said {\em of rank one} (resp.
{\em binomial}) if $\sigma_t(x)=(1-tx)^{-1}$ (resp.
$\sigma_t(x)=(1-t)^{-x}$). The scalar $x=1$ is the only element
having both properties. We usually consider that our auxiliary variable
$t$ is of rank one, so that $\sigma_t(A)=\sigma_1(tA)$.

The argument $A$ of the noncommutative symmetric functions can be a
``virtual alphabet''. This means that, being algebraically independent,
the $S_n$ can be specialized to any sequence $\alpha_n\in{\mathcal A}$ of
elements of any associative algebra ${\mathcal A}$. Writing $\alpha_n=S_n(A)$
defines all the symmetric functions of $A$, and allows one to
use the powerful notations $F(nA)$, $F((1-q)A)$, etc., for more
or less complicated transformations of the specialized functions.

The (commutative) specializations $A=\E$, defined by
\begin{equation}
S_n(\E)=\frac{1}{n!}
\end{equation}
and $A=\frac{1}{1-q}$, for which
\begin{equation}
S_n\left(\frac{1}{1-q}\right)=\frac{1}{(q)_n}
\end{equation}
are of special importance.

\subsection{Noncommutative analogs of special functions}

Since the discovery by D. Andr\'e of a combinatorial interpretation of tangent
and secant numbers, several classical generating functions have been lifted
to the algebra of symmetric functions, and more recently, to noncommutative
symmetric functions. The general idea is as follows. Given the exponential
generating function 
\begin{equation}
f(t)=\sum_{n\ge 0}c_n \frac{t^n}{n!}
\end{equation}
of a combinatorial sequence $a_n\in\N$, one looks for a noncommutative
symmetric function $F(A)$ such that $F(t\E)=f(t)$. The noncommutative
analog is interesting when $F_n(A)$ can be directly interpreted as the formal
sum of the combinatorial objects counted by  $c_n$, under the embedding
of $\Sym$ into some larger algebra. For example, in the case of tangent
and secant
numbers, the series
\begin{equation}
\left(\sum_{n\ge 0}(-1)^nS_{2n}(A)\right)^{-1}
\left(1+ \sum_{n\ge 0}(-1)^nS_{2n+1}(A)\right)
\end{equation}
becomes the formal sum of the alternating permutations 
(shapes $(2^n)$ and $(2^n1)$) under the embedding of $\Sym$
in $\FQSym$ \cite{NCSF1}. One can also find in \cite{NCSF1}
the noncommutative Eulerian polynomials, and in \cite{NT3},
analogs of the Abel polynomials and of the Lambert and Eisenstein
functions.

In general,  $F_n$ turns out to be the characteristic
of some projective $0$-Hecke module. Projective
modules are always specializations of generic modules, thus also
representations of the symmetric group, whose Frobenius characteristic
are then the commutative images $F_n(X)$. In general, setting
$X=\frac{t}{1-q}$ gives back an interesting $q$-analog of $f(t)$.

In this note, we shall show that the consideration of $0$-Hecke
modules obtained from a natural notion of tensor products leads
immediately to noncommutative anlogs of the Bessel $J$ (or $I$)
functions. Here, we need two alphabets $A$ and $B$, and we are
led to the combinaorics of bi-exponential generating functions.

\section{Tensor products of $0$-Hecke modules}

\subsection{}
The $0$-Hecke algebra $H_n(0)$ is the  algebra
$\C[\Pi_n]$, where the monoid $\Pi_n$ is generated by elements
$\pi_1,\ldots,\pi_{n-1}$ ($\pi_i=1+T_i$) satifying the braid relations
\begin{eqnarray}
\pi_i\pi_j&= \pi_j\pi_i & |i-j|>1\\
\pi_i\pi_{i+1}\pi_i &=\pi_{i+1}\pi_i\pi_{i+1}&\\
\end{eqnarray}
and the idempotency condition
\begin{equation}
\pi_i^2=\pi_i \,.
\end{equation}
There is a canonical coproduct on $H_n(0)$ defined by
\begin{equation}
\delta_\wedge\pi= \pi\otimes\pi \ \text{for $\pi\in\Pi_n$}\,.
\end{equation}
Hence, tensor products of $H_n(0)$-modules can be defined,
and it is obvious from the definition of the simple module
$S_I$ that
\begin{equation}
S_H\otimes S_K =S_I \quad\text{where $\Des(I)=\Des(H)\cap\Des(K)$}\,.
\end{equation}
This induces an internal product $\wedge$ on $QSym_n=G_0(H_n(0))$, similar
to the internal product of symmetric functions, such that
\begin{equation}
F_H \wedge F_K =F_I \quad\text{where $I=H\wedge K$, that is, $\Des(I)=\Des(H)\cap\Des(K)$}\,.
\end{equation}
By duality, this defines a coproduct on $\Sym_n$, given by
\begin{equation}
\gamma_\wedge R_I = \sum_{\Des(I)=\Des(H)\cap\Des(K)}R_H\otimes R_K\,.
\end{equation}

\subsection{} 
There is a canonical involution $\iota$ on $H_n(0)$, defined by 
\begin{equation}
\iota(\pi)=\bar\pi_i = 1-\pi_i\,,
\end{equation}
so that we can regard $H_n(0)$ as $\C[\bar\Pi_n]$ as well.
Hence, we have another tensor product, defined from the coproduct
\begin{equation}
\delta_\vee\bar\pi= \bar\pi\otimes\bar\pi \ \text{for $\pi\in\Pi_n$}\,,
\end{equation}
which induces a second internal product $\vee$ on $QSym$,
\begin{equation}
F_H \vee F_K =F_I \quad\text{where $I=H\vee K$, that is,
$\Des(I)=\Des(H)\cup\Des(K)$}\,.
\end{equation}
It is of course sufficient to study one of them. However,
it is interesting to observe that this second product
appears in another guise in \cite{NCSF6}, in the process
of calculating a basis of primitive elements of $\FQSym$.
Let us recall this construction. 
Let $p_n$ denote the projection onto the homogeneous
component $\FQSym_n$ of $\FQSym$, and let $\mu_q:\
\F_\alpha\otimes\F_\beta \mapsto \F_{\alpha\shuffle_q \beta[k]}$ be
the multiplication map of $\FQSym_q$. The $q$-convolution of two
graded linear endomorphisms $f,g$ of $\FQSym$ is defined by
\begin{equation}
f \odot_q g = \mu_q\circ(f\otimes g)\circ\Delta\, .
\end{equation}
For $q=1$, this reduces to ordinary convolution.
We are interested
in the case $q=0$. For a composition $I=(i_1,\ldots,i_m)$, let
\begin{equation}
p_I= p_{i_1}\odot_0\cdots\odot_0 p_{i_m} \,.
\end{equation}
It is proved in \cite{NCSF6} that
the $p_I$ are mutually commuting projectors, and more precisely that
  \begin{equation*}
    p_I\circ p_J= \left\{
      \begin{array}{cl}
        0  & \ \hbox{if} \ |I|\not= |J|. \\[1mm]
        p_{I\vee J} & \ \hbox{otherwise} \ .
      \end{array}
    \right.
\end{equation*}
Hence, $j:\ F_I\mapsto p_I$ defines an embedding of $(QSym,\vee)$
in the composition algebra of graded endomorphisms of $\FQSym$.
Moreover, 
\begin{equation}
\pi=\sum_{|I|\ge 1} (-1)^{l(I)-1} p_I
\end{equation}
which is a projector onto the primitive Lie algebra of $\FQSym$,
is the image of the primitive element $\sum_nM_n$ of $QSym$ under $j$,
and it easy to see that more generally, for any $f\in QSym$ 
\begin{equation}
(j\otimes j)(\Delta_{QSym}f) = \Delta_{\FQSym}\circ j(f)\,.
\end{equation}
However, $j$ does not map the usual (external) product of $QSym$
to the ordinary convolution of endomorphisms. It is nevertheless
interesting to pull back the $0$-convolution to $QSym$, by
defining
\begin{equation}
F_I \odot_0 F_J = F_{I\cdot J}\,,
\end{equation}
where $I\cdot J$ means as usual concatenation of the compositions.
Then, we have a splitting formula
\begin{equation}
(f_1\odot_0 f_2 \odot_0\cdots\odot_0 f_r)\vee g
=
\mu_0[(f_1\otimes\cdots\otimes f_r)\vee_r \Delta_{QSym}^r(g)]        
\end{equation}
analogous to the one satisfied in $\Sym$.

It can be shown that the involution $\iota$ maps the simple
module $S_I$ and the indecomposable projective module $P_I$
to $S_{\bar I^\sim}$ and $P_{\bar I^\sim}$, respectively.

\subsection{}
Identifying as usual a tensor product $F\otimes G$ with $F(A)G(B)$,
where $A$ and $B$ are two mutually commuting alphabets, we have
\begin{equation}
\sigma_1(XA)\wedge\sigma_1(XB) = \sum_K F_K(X)\gamma_\wedge(R_K)=\gamma_\wedge\sigma_1(XA)\,,
\end{equation}
which may be compared with the following identity relating the internal
product $*$ of $\Sym$ and its dual coproduct $\delta F=F(XY)$ on $QSym$:
\begin{equation}
\sigma_1(XA)*\sigma_1(YA)=\sigma_1(XYA)=\delta\sigma_1(XA)\,.
\end{equation}

\begin{theorem}
The coproduct $\gamma_\wedge$ is a morphism for the ordinary (outer) product
of non commutative symmetric functions, that is
\begin{equation}
\gamma_\wedge(FG)=\gamma_\wedge(F)\gamma_\wedge(G)
\end{equation}
In particular, it is completely determined by the images of the elementary
functions, $\gamma_\wedge\Lambda_n=\Lambda_n\otimes\Lambda_n$, which implies the
combinatorial inversion formula
\begin{equation}
\left(\sum_{n\ge 0}(-1)^n\Lambda_n\otimes\Lambda_n\right)^{-1}
=\sum_{\Des(H)\cap\Des(K)=\emptyset}R_H\otimes R_K\,.
\end{equation}
\end{theorem}

\Proof This is equivalent to Theorem \ref{thRub} below. \qed

As we will see, this simple identity has many interesting enumerative
corollaries. Applying the involution $\omega$ on the second factor
gives the inverse of
\begin{equation}
\left(\sum_{n\ge 0}(-1)^n\Lambda_n\otimes S_n\right)^{-1}
=\sum_{\Des(H)\cap\Des(K)=\emptyset}R_H\otimes R_{K^\sim}\,.
\end{equation}
The right hand side of this equality occurs in \cite{HT},
where it is interpreted as the decomposition of the 
algebra ${\rm H}\SG_n$ as a bimodule over itself.
The inverse of the left hand side  can legitimately be considered
as a noncommutative analog of the Bessel function $J_0$, as if
we specialize both sides to $x\E$, we recover $J_0(2x)$.
Moreover, specializing $A$ to $x/(1-q)$ gives a classical $q$-analog
of $J_0$, and the other ones are obtained by simple transformations.
This first step being granted, it is not difficult to guess the
correct definition of the noncommutative analogues of the other
$J_\nu$. This will be done in the forthcoming section.

\section{Noncommutative Bessel functions}

\subsection{} Let $A$ and $B$ be two mutually commuting alphabets.
The noncommutative Bessel functions $\J_n(A,B)$ are defined by the
generating series
\begin{equation}
\sum_{n\in\Z}z^n \J_n(A,B)=\lambda_{-1/z}(A)\sigma_z(B)\\,
\end{equation}
that is,
\begin{equation}
\J_n(A,B)=\sum_{m\ge 0}(-1)^m\Lambda_{m-n}(A)S_m(B)\,.
\end{equation}
For $A=B=x\E$, this is the usual Bessel function $J_n(2x)$.
In particular,
\begin{equation}
\J_0(A,B)=\sum_{m\ge 0}(-1)^m\Lambda_{m}(A)S_m(B)
\end{equation}
can be regarded as $\lambda_{-1}(\JJ)$, for  the virtual alphabet
$\JJ=(A,B)$ such that
\begin{equation}
\Lambda_n(\JJ)=\Lambda_n(A)S_n(B)\,.
\end{equation}
This defines an embedding of algebras
\begin{equation}
\begin{split}
j:\ \Sym & \rightarrow \Sym(A,B)=\Sym\otimes\Sym\\
\Lambda_n(A)&\mapsto \Lambda_n(\JJ)=\Lambda_n\otimes S_n\,.
\end{split}
\end{equation}
It is not difficult to describe the image of the ribbon basis
under this embedding. We need the following piece of notation.
For two compositions $I$ and $J$ of the same integer $n$,
we define the composition $K=I\backslash J$ of $n$ by the
condition
\begin{equation}
\Des(K)=\Des(I)\backslash \Des(J)\quad\text{(set difference)}\,.
\end{equation}
Then, we can state:
\begin{theorem}\label{thRub}
The image of $R_K$ by $j$ is
\begin{equation}
R_K(\JJ)=\sum_{I\backslash J=K}R_I(A)R_J(B)\,.
\end{equation}
\end{theorem}

\Proof The formula is true for $K=(1^n)$ by definition. The
general case follows by induction on $l(K^\sim)$, the number of columns
of the ribbon diagram of $K$. Indeed, it suffices to prove that
\begin{equation}
R_K(\JJ)R_{1^m}(\JJ)=R_{K1 ^m}(\JJ)+R_{K\triangleright 1^m}(\JJ)\,,
\end{equation}
which follows easily from the usual multiplication rule of ribbon functions.
\qed

\begin{corollary}[\cite{CSV}]
Let $a_n$ be defined by
\begin{equation}
\frac{1}{J_0(2\sqrt{t})}=\sum_{n\ge 0}a_n\frac{t^n}{(n!)^2}\,.
\end{equation}
Then, $a_n$ is equal to the number of pairs of permutations
$(\sigma,\tau)\in\SG_n\times\SG_n$ such that $\Des(\sigma)\subseteq
\Des(\tau)$.
\end{corollary}

Let $\dd$ be the linear operator on $\Sym$ (acting on the right)
defined by
\begin{equation}\label{defdd}
S^{(i_1,\ldots,i_r)}\dd=S^{(i_1,\ldots,i_r-1)}\,.
\end{equation}
It has the following properties (see \cite{NCSF2}, Prop. 9.1):
\begin{equation}
(FG)\dd = F\cdot (G\dd)+ (F\dd)\cdot G_0\,,
\end{equation}
where $G_0$ denotes the constant term of $G$, and
\begin{equation}
R_I\dd=\left\{\begin{matrix}R_{i_1,\ldots,i_r-1}&\text{if $i_r>1$}\,,\\
0&\text{if $i_r=1$}\end{matrix}\right.\,.
\end{equation}
In particular, if $G_0=0$,
\begin{equation}
(1-G)^{-1}\dd = (1-G)^{-1} (G\dd)\,.
\end{equation}
Let us apply this with $\dd=\dd_B$ acting only on $\Sym(B)$ to
\begin{equation}
\J_0(A,B)^{-1}=
\left( 1-\sum_{n\ge 1}(-1)^{n-1}\Lambda_n(B)S_n(A)\right)^{-1}
=\sum_I S^I(A)R_I(B)\,.
\end{equation}
We obtain
\begin{equation}
\J_0(A,B)^{-1}\J_{-1}(A,B)=\sum_I S^I(A)(R_I\dd(B))
\end{equation}

\begin{corollary}[\cite{CSV}]
The coefficient $c_n$ in
\begin{equation}
\frac{J_1(2{x})}{J_0(2{x})}=\sum_{n\ge 1}c_n \frac{x^{2n-1}}{(n-1)!n!}
\end{equation}
is equal to the number of pairs of permutations $(\alpha,\beta)\in\SG_n^2$
such that $\Des(\alpha)\subseteq\Des(\beta)$ and $\beta(n)=n$.
\end{corollary}


\subsection{Bessel-Carlitz functors}

Let $\F$ be the functor which associates with a pair
of vector spaces $(V,W)$ the graded subalgebra of
the exterior algebra $\Lambda(V\oplus W)$
\begin{equation}
\F(V,W)=\bigoplus_{n\ge 0} \Lambda_n(V)\otimes \Lambda_n(W)\,.
\end{equation}
This is a quadratic algebra (see\cite{Pr}). 
If $(v_i)$, $(w_j)$ are bases of $V$ and of
$W$, the relations are as follows. For $i<k$ and $j<l$,
\begin{eqnarray} 
\biw{i\,k}{j\l}+\biw{k\,i}{j\,l}& = 0\,,\\
\biw{i\,k}{j\l}+\biw{i\,k}{l\,j}& = 0\,,\\
\biw{i\,i}{j\,l}& = 0\,,\\
\biw{i\,k}{j\,j}& = 0\,,
\end{eqnarray}
where $\biw{i}{k}=v_i\otimes w_k$.

Hence,  the Koszul dual $\G(V,W)=\F(V,W)^!$ is the quadratic algebra
on $V^*\otimes W^*$ presented by 
\begin{equation}
\biw{i\,k}{j\l}=\biw{k\,i}{j\,l}=\biw{i\,k}{l\,j}\ \text{for $i<k$ and $j<l$.}
\end{equation}

The combinatorial investigation of Bessel functions has
been initiated by Carlitz \cite{Ca}. Hence,
the polynomial bi-functors defined by $\F$ and $\G$ can appropriately
be called Bessel-Carlitz functors. One or two occurences of $\Lambda$
can be replaced by $S$ in the definition of
$\F$. In the mixed case $\Lambda\otimes S$, the
best interpretation is probably as functors defined on super (i.e.,
$\Z_2$-graded) vector spaces $V=V_0\oplus V_1$.


\section{The $\theta$-specialization}

This section is devoted to the interpretation of a few formulas
from \cite{CSV2,FR2,FR3} in terms of noncommutative symmetric functions.

\subsection{}Let $\theta\subseteq A\times A$ be any binary relation.
We denote by $\overline\theta$ the complement of $\theta$ in
$A\times A$ and set
\begin{equation}
\begin{split}
X=&X(A;\theta)=\{w=a_1\cdots a_n\in A^*|a_1 \theta a_2\theta\ldots\theta
a_n\}\,,\\
Y=&Y(A;\theta)=X(A;\overline\theta)\,,
\end{split}
\end{equation}
where we write $a\theta b$ for $(a,b)\in\theta$.
Note that the empty word $1$ and the letters belong to both $X$ and $Y$.

The $\theta$-specialization $\Sym(A;\theta)$ is then defined by
specifying the elementary symmetric functions
\begin{equation}
\Lambda_n(A;\theta)=\sum_{w\in X\cap A^n} w\,.
\end{equation}
The following basic lemma,
implicit in \cite{CSV2},
 generalizes the case $\theta=\{(a,b)|a>b\}$.
\begin{lemma}[Carlitz-Koszul duality for alphabets]
The complete symmetric functions $S_n(A;\theta)$ are given by
\begin{equation}
S_n(A;\theta)=\Lambda_n(A;\overline\theta)\,.
\end{equation} 
More generally, if one denotes by $\tAdj(w)=\{i|a_i\theta a_{i+1}\}$
the $\theta$-adjacency set of $w=a_1a_2\cdots a_n$, and by $C_\theta(w)$
the associated composition of $n$, one has
\begin{equation}
R_I(A;\theta)=\sum_{C_\theta(w)=I}w\,.
\end{equation} 
\end{lemma}

\Proof 
We need to prove that
\begin{equation}
\sum_{k=0}^n(-1)^k\Lambda_k(A,\theta)\Lambda_{n-k}(A,\bar\theta)=0
\end{equation}
for $n>0$. Let $w=uv$ be such that $u\in\Lambda_k(A,\theta)$ and 
$v\in \Lambda_{n-k}(A,\bar\theta)$. Then if ${\rm last}(u)\theta
{\rm first}(v)$, $w$ appears in
$\Lambda_{k+1}(A,\theta)\Lambda_{n-k-1}(A,\bar\theta)$,
and similarly, if ${\rm last}(u)\bar\theta
{\rm first}(v)$, then $w$ appears in
$\Lambda_{k-1}(A,\theta)\Lambda_{n-k+1}(A,\bar\theta)$.
Moreover, $w$ cannot appear in any other product, so that its coefficient
in the sum is $0$. \qed

\subsection{The $\theta$-Eulerian polynomials}
Recall from \cite{NCSF1} that the noncommutative Eulerian polynomials
\begin{equation}
\A_n(t;A)=\sum_{I\vDash n}t^{l(I)}R_I(A)
\end{equation}
admit the generating function
\begin{equation}
\AA(t;A)=\sum_{n\ge 0}\A_n(t;A)=\frac{1-t}{1-t\sigma_{1-t}(A)}
\end{equation}
(see \cite{Desar} for the commutative version of this identity),
and since $l(C_\theta(w))=\tadj(w)+1$, we have immediately
\begin{equation}
\sum_{w\in A^*}t^{\tadj(w)+1}w=\frac{1-t}{1-t\sigma_{1-t}(A;\theta)}\,.
\end{equation}
Note that $\tadj(w)+{\overline{\theta}}{\rm adj}(w)=n-1$.
Replacing $\theta$ by ${\overline\theta}$, $A$ by $t^{-1}A$, then $t$
by $t^{-1}$, and simplifying by $(1-t)$ the resulting expression,
we obtain
\begin{equation}\label{th2deFR}
\sum_{w\in A^*}t^{\tadj(w)}w = 
\frac{1}{1-\sum_{w\in X(A;\theta)^{+}} (t-1)^{l(w)-1} w}\,,
\end{equation}
which is Theorem 2 of \cite{FR3}. 

For a letter $c\in A$, denote by $\dd_c$ the linear operator
defined by
\begin{equation}
w\dd_c = \left\{\begin{matrix} u & \text{if $w=uc$ for some $u$}\,,\\
                                    0 & \text{otherwise}\,.
\end{matrix}\right.
\end{equation}
Then, as in (\ref{defdd}), for any series $F$ without constant term,
\begin{equation}
(1-F)^{-1}\dd_c = (1-F)^{-1}\cdot (F\dd_c)\,.
\end{equation}
The same is true for the operators
\begin{equation}
D_C=\sum_{c\in C} \dd_c \cdot c
\end{equation}
where $C$ is a subset of $A$. Applying this to (\ref{th2deFR}),
we obtain
\begin{equation}\label{th3deFR}
\sum_{w\in A^*C}t^{\tadj(w)}w = \frac{-\sum_{w\in XC}(t-1)^{l(w)-1}w}
{1-\sum_{w\in X(A;\theta)^+}
  (t-1)^{l(w)-1}w}\,,
\end{equation}
which is Theorem 3 of \cite{FR3}.

\subsection{The $\theta$-Major index}
If one defines the $\theta$-Major index by
\begin{equation}
\tmaj(w)=\sum_{i\in\tAdj(w)}i
\end{equation}
one has clearly
\begin{equation}
\sum_{w\in A^n}q^{\tmaj(w)}w=\sum_{I\vdash n}q^{\maj(I)}R_I(A;\theta)
=(q)_nS_n\left(\frac{A}{1-q};\theta\right)\,,
\end{equation}
where as usual
\begin{equation}
\sigma_z\left(\frac{A}{1-q};\theta\right)=\prod_{n\ge
0}^{\rightarrow}\sigma_{zq^n}(A;\theta)\,.
\end{equation}

\section{Double Eulerian polynomials and Bessel functions}

\subsection{}The noncommutative Bessel function $\J_0(A,B)$
can now be properly interpreted as a generating series of 
$\theta$ elementary symmetric functions, if we interpret
$\JJ$ as the product alphabet $A\times B$, endowed with
the relation
\begin{equation} 
(a,b)\theta (a',b')\ \Leftrightarrow \ a>a'\ \text{and}\ b\le b'\,.
\end{equation}
As is customary, we denote words over $A\times B$ by biwords
\begin{equation}
w=[u,v]=\biw{u}{v}\,\ u\in A^n\,,\ v\in B^n\,.
\end{equation}
Observing that
\begin{equation}
\tAdj\left(\biw{u}{v}\right) = \Des(u)\cap
\overline{\Des(v)}=\Des(u)\backslash\Des(v)\,,
\end{equation}
we can now write
\begin{equation}\label{doubleEuler}
\begin{split}
\sum_{w=(u,v)\in(A\times B)^*}t^{\tadj(w)}z^{l(w)} w
&=\frac{1-t}{\J_0((1-t)z;A,B)-t}\\
&=\sum_K z^{|K|}t^{l(K)-1}R_K(A,B;\theta)
\end{split}
\end{equation}
where from now on we shall use the notation
\begin{equation}
\J_0(x;A,B)=\lambda_{-x}(\JJ)=\lambda_{-x}(A,B;\theta)\,.
\end{equation}
The coefficient of $z^n$ is the $n$th double $\theta$-Eulerian
polynomial, denoted by $\A_n(t;A,B;\theta)$. Setting $A=B=\E$,
we recover the enumeration of pairs of permutations
$(\alpha,\beta)\in\SG_n\times\SG_n$ by the cardinality of
$\Des(\alpha)\cap \overline{\Des(\beta)}$ (cf. \cite{CSV}).

\section{The F\'edou-Rawlings polynomials}
By considering simultaneously
the specializations of (\ref{doubleEuler})
to all positive $q$ and $p$-integers, $A_i=[i+1]_q$ and $B_j=[j+1]_p$,
one arrives at the five parameter generalizations of the double Eulerian
polynomials introduced by F\'edou and Rawlings \cite{FR3}.

For $w\in A^n$, where $A$ is the infinite
chain $A=\{a_1<a_2<\ldots\}$, let $q^w$ be the image
of $w$ by the multiplicative homomorphism $a_i\mapsto q^{i-1}$.
Writing, for a composition $I$ of $n$
\begin{equation}
R_I(A) = \sum_{C(\sigma)=I} \sum_{\Std(w)=\sigma}w
\end{equation}
and taking into account the  identity 
\begin{equation}
\sum_{\Std(w)=\sigma}x^{\max(w)}q^w = \frac{x^{\des(\sigma^{-1})}q^{\coimaj(\sigma)}}{(xq;q)_n}
\end{equation}
where $\coimaj(\sigma)$ denotes the co-major index of $\sigma^{-1}$,
\begin{equation}
\coimaj(\sigma)=\sum_{d\in\Des(\sigma^{-1})}(n-d)\,,
\end{equation}
(indeed, it is easily checked that the minimal word $v$ for the
lexicographic order such that $\Std(v)=\sigma$ satisfies
$q^v=q^{\coimaj(\sigma)}$), we find
\begin{equation}
\sum_{i\ge 0} x^i R_I(1,q,\ldots,q^i)=\frac{1}{1-x}\sum_{C(w)=I}x^{\max(w)}q^w
=\frac{1}{(x;q)_{n+1}}\sum_{C(\sigma)=I}x^{\des(\sigma^{-1})}q^{\coimaj(\sigma)}
\end{equation}
so that finally,
we recover the double generating series of \cite{FR3}
\begin{multline}\label{FRserie1}
\sum_{i,j\ge 0}x^iy^j \frac{1-t}{\J_0((1-t)z;A_i,B_j)-t}\\
=
\sum_{n\ge 0}\frac{z^n}{ (x;q)_{n+1} (y;p)_{n+1} }
\sum_{\alpha,\beta\in\SG_n}
 { t^{\desris(\alpha,\beta)}
x^{\des(\alpha^{-1})}y^{\des(\beta^{-1})}
q^{\coimaj(\alpha)} p^{\coimaj(\beta)} } \,,
\end{multline}
where 
$\desris(\alpha,\beta)=|\Des(\alpha)\backslash\Des(\beta)|$.

The second generating series of \cite{FR3} is recovered in the same 
way. If we denote by $b_j$ the greatest letter of $B_j$, then,
on the one hand,
\begin{equation}
S_n(B_j)\dd_{b_j} = S_{n-1}(B_j)\,.
\end{equation}
On the other hand,
\begin{multline}
\sum_{j\ge 0}
y^j R_J(B_j)\,\dd_{b_j}\cdot b_j=
\frac{1}{1-y}
\sum_{\empile{C(\sigma)=J}{ \max(v)=\last(v)}}y^{max(v) }v\\
=\frac{1}{1-y}\sum_{\empile{C(\sigma)=J}{ \sigma(n)=n}} 
\sum_{\Std(v)=\sigma}y^{max(v) }v\,,
\end{multline}
so that, applying the operator $\dd_{b_j}\cdot b_j$
to the coefficient of $x^iy^j$ in (\ref{FRserie1}), we obtain
\begin{multline}
\sum_Kz^{|K|}\sum_{I\backslash J=K} t^{l(K)-1}\sum_{i,j\ge 0}x^iy^j
R_I(A_i)R_J(B_j)\cdot\dd_{b_j}\cdot b_j\\
=\sum_{i,j\ge 0}x^iy^j
\left(
1-\sum_{n\ge 1}z^n(t-1)^{n-1}\Lambda_n(A_i)S_n(B_j)
\right)^{-1}\cdot\dd_{b_j}\cdot b_j\\
=\sum_{i,j\ge 0}x^iy^j
\frac{
\left(
-\sum_{n\ge 1}z^n(t-1)^{n-1}\Lambda_n(A_i)S_{n-1}(B_j)b_j
\right)(1-t)}
{\J_0(z(1-t);A_i,B_j)-t}\\
=\sum_{i,j\ge 0}x^iy^j
\frac{
\J_{-1}((1-t)z;A_i,B_j)b_j}
{\J_0(z(1-t);A_i,B_j)-t}\,.
\end{multline}
Specializing $A_i=[i+1]_q$, $B_j=[j+1]_p$, this becomes,
in the notation of \cite{FR3},
\begin{multline}
\sum_{i,j\ge 0}x^i(py)^j
\frac{J_1^{(i,j)}((1-t)z;q,p)}{J_0^{(i,j)}((1-t)z;q,p)-t}\\
=
\sum_{n\ge 0}\frac{z^n}{(x;q)_{n+1}(y;p)_n}
\sum_{\empile{\alpha,\beta\in\SG_n}{ \beta(n)=n}}
t^{\desris(\alpha,\beta)}
x^{\des(\alpha^{-1})}
y^{\des(\beta^{-1})}
q^{\coimaj(\alpha)}
p^{\coimaj(\beta)}
\end{multline}
which is equivalent to \cite[(3)]{FR3}. 
Here,
\begin{equation}
J_\nu^{(i,j)}(z;q,p):=(-1)^\nu\J_\nu(z[i+1]_q,[j+1]_p)\,.
\end{equation}
The other results
of \cite{FR3} can be rederived in the same way, by changing
the specializations of $A_i$ and $B_j$.

\section{Heaps of segments and polyominos}

Bessel functions and their multiparameter analogs play a crucial role
in the enumerative theory of polyominos \cite{BMV,DF}. 
Elegant combinatorial proofs
of such enumerative results can be achieved by means of Viennot's
theory of heaps of segments \cite{Vi,BMV}. As we shall see,
this can also be  conveniently formulated in terms of
$\theta$-noncommutative symmetric functions.

\subsection{}
A {\em parallelogram} (or {\em staircase}) {\em polyomino} $P$,
which is also the same as a connected skew Young diagram, can be encoded as a
biword
\begin{equation}
w=a_{i_1j_1}\cdots a_{i_nj_n}=\biw{i_1\cdots i_n}{j_1\cdots j_n}=\biw{u}{v}
\end{equation}
where $j_k$ is the heigth of the $k$th column $C_k$, and $i_k$ is the
number of common rows between $C_k$ and $C_{k+1}$ (with a conventional
value $i_n=1$ for the last column). For example, the following polyomino

\centerline{\includegraphics[height=5cm]{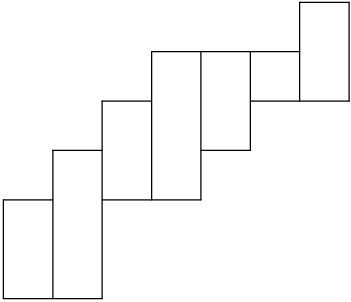}}
is encoded by the biword
\begin{equation}
\biw{ 2 1 2 2 1 1 1}{ 2 3 2 3 2 1 2}
\end{equation}
The biwords corresponding to polyominos are the words over the alphabet
\begin{equation}
A=\{a_{ij}|i\le j\}
\end{equation}
satisfying the $\theta$-adjacency conditions
\begin{equation}\label{adjpol}
a_{i_k j_k}\,\theta\, a_{i_{k+1} j_{k+1}}\ \Longleftrightarrow i_k \le j_{k+1}
\end{equation}
and the ending condition
\begin{equation}\label{endpol}
i_n=1\,.
\end{equation}
Hence, the generating series (by length)  of all biwords satisfying (\ref{adjpol}) is
\begin{equation}
\lambda_t(A,\theta)=[\lambda_{-t}(A,\bar\theta)]^{-1}
=\left(1-\sum_{n\ge 1}(-1)^{n-1}t^n \sum_{i_k>j_{k+1}}  \biw{i_1\cdots i_n}{j_1\cdots j_n}\right)^{-1}
\end{equation}
and restriction of the series to the biwords satifying (\ref{endpol}) is
achieved as above by applying the operator
\begin{equation}
D=\sum_{j\ge 1}\dd_{a_{1j}}
\end{equation}
so that we end up once more with a series of the form
\begin{equation}
\left(1-\sum_{n\ge 1}(-1)^{n-1}t^n \sum_{i_k>j_{k+1}}  \biw{i_1\cdots
i_n}{j_1\cdots j_n}\right)^{-1}
\left(-\sum_{n\ge 1}(-1)^{n-1}t^n \sum_{i_k>j_{k+1};i_n=1}  \biw{i_1\cdots
i_n}{j_1\cdots j_n}\right)
\end{equation}
which acquires the structure $J_1/J_0$ once $A$ is specialized to 
\begin{equation}
a_{ij}=x y^{j-i}q^j\,,
\end{equation}
the generating series by  width, height and area.

\subsection{} This can of course be interpreted in terms of heaps of segments.
A {\em segment} is an interval $[i,j]$ of $\N^*$. To each segment,
we associate a variable
\begin{equation}
a_{ij}=\biw{i}{j}\,,
\end{equation}
in our $A=\{a_{ij}|i\le j\}$. The {\em monoid of heaps} is the quotient
of the free monoid $A^*$ by the commutation relations
\begin{equation}
a_{ij}a_{kl}\equiv a_{kl}a_{ij}\ \text{if $j<k$}
\end{equation}
which means that the segments do not overlap and can be vertically
slided independently of each other.

The first basic lemma of the theory (which is also a special case of 
the Cartier-Foata formula for the Moebius functions of free partially
commutative monoids \cite{CF}) amounts to the calculation of $S_n(A,\theta)$
for the relation defined by
\begin{equation}
a_{ij}\,\theta\,a_{kl} \Longleftrightarrow i\le l\,.
\end{equation}
Indeed, with this choice, $\Lambda_n(A,\bar\theta)$ is the formal sum of trivial
heaps (products of mutually commuting segments arranged in decreasing order),
and $\Lambda_n(A,\theta)$ is the sum of all biwords
\begin{equation}
w=a_{i_1j_1}\cdots a_{i_nj_n}=\biw{i_1\cdots i_n}{j_1\cdots j_n}
\end{equation}
such that $i_k\le j_{k+1}$ for all $k$, those encoding polyominos. 
We have therefore shown that each heap, or, equivalently, each element of $A^*/\equiv$ has a unique
representative of this form.


\end{document}